\documentclass[11pt]{amsart}
\usepackage{thm-restate,amsmath,amssymb,enumerate,graphicx}
\usepackage[
  colorlinks,
  linkcolor = blue,
  citecolor = blue,
  urlcolor = blue]{hyperref}
\usepackage{tikz}
\usepackage{subfig}
\newtheorem{remark}{Remark}[section]
\newtheoremstyle{plainsl}%
	{\topsep}
	{\topsep}
	{\slshape} 
	{}
	{\normalfont\bfseries}
	{.}
	{ }
	{}
\swapnumbers

\newtheorem{theorem}{Theorem}[section]
\newtheorem{lemma}[theorem]{Lemma}
\newtheorem{corollary}[theorem]{Corollary}
\newtheorem{conjecture}[theorem]{Conjecture}

\newtheorem*{cor1}{Corollary}
\newtheorem*{thm2}{Theorem}
\newtheorem*{thm3}{Theorem}

\newcommand\cref[1]{Corollary~\ref{cor:#1}}

\renewcommand\proof{\noindent\textsl{Proof. }}
\newcommand\sqr[2]{{\vbox{\hrule height.#2pt
    \hbox{\vrule width.#2pt height#1pt \kern#1pt
        \vrule width.#2pt}\hrule height.#2pt}}}
\renewcommand\qed{%
	\ifmmode\eqno\sqr53
	\else\nolinebreak\ \hfill\sqr53\medbreak\fi}

\begin{document}

\sloppy
\title[Entropy of Symmetric Graphs]{Entropy of Symmetric Graphs}
%
%
%


\author[Changiz Rezaei, and Godsil]{Seyed Saeed Changiz Rezaei, and Chris Godsil}

\address{
        Seyed Saeed Changiz Rezaei, Chris Godsil,
        Department of Combinatorics and Optimization,
        University of Waterloo,
        Waterloo, Canada}
\email{\{sschangi,cgodsil\}@uwaterloo.math.ca}

\date{\today}
\maketitle

\begin{abstract}
 A graph $G$ is called \emph{symmetric with respect to a functional $F_G(P)$} defined on the set of all the probability distributions on its vertex set if the distribution $P^*$ maximizing $F_G(P)$ is uniform on $V(G)$. Using the combinatorial definition of the entropy of a graph in terms of its vertex packing polytope and the relationship between the graph entropy and fractional chromatic number, we prove that vertex transitive graphs are symmetric with respect to graph entropy. As the main result of this paper, we prove that a perfect graph is symmetric with respect to graph entropy if and only if its vertices can be covered by disjoint copies of its maximum-size clique. Particularly, this means that a bipartite graph is symmetric with respect to graph entropy if and only if it has a perfect matching.
\end{abstract}

\section{Introduction}
The entropy of a graph is an information theoretic functional which is defined on a graph with a probability density on its vertex set. This functional was originally proposed by J. K\"{o}rner in 1973 to study the minimum number of codewords required for representing an information source (see J. K\"{o}rner \cite{JKor}).

Let $VP(G)$ be the \emph{vertex packing polytope} of a given graph $G$ which is the convex hull of the characteristic vectors of its independent sets. Let $|V(G)| = n$ and $P$ be a probability density on $V(G)$. Then the \emph{entropy of $G$ with respect to the probability density $P$} is defined as
\[
H(G,P) = \min_{\mathbf{a}\in VP(G)} \sum_{i=1}^n p_i\log (1/a_i).\label{eq:combent}
\]

G. Simonyi \cite{Simu} showed that the maximum of the graph entropy of a given graph over the probability density of its vertex set is equal to its fractional chromatic number. We say a graph is \emph{symmetric with respect to graph entropy} if the uniform density maximizes its entropy. We show that vertex transitive graphs are symmetric. In this paper, we study some classes of graphs which are symmetric with respect to graph entropy.  Our main results are the following theorems and corollary.

\begin{thm2}
Let $G=(V,E)$ be a perfect graph and $P$ be a probability distribution on $V(G)$. Then $G$ is symmetric with respect to graph entropy $H\left(G,P\right)$ if and only if $G$ can be covered by disjoint copies of its maximum-size cliques.
\end{thm2}
As a corollary to above theorem, we have
\begin{cor1}
Let $G$ be a bipartite graph with parts $A$ and $B$, and no isolated vertices. Then, uniform probability distribution $U$ over the vertices of $G$ maximizes $H\left(G,P\right)$ if and only if $G$ has a perfect matching.
\end{cor1}
A. Schrijver \cite{Schriv1} calls a graph $G$ a \emph{$k$-graph}\index{$k$-graph} if it is $k$-regular and its fractional edge coloring number $\chi_f^\prime(G)$ is equal to $k$. We show that
\begin{thm3}
Let $G$ be a $k$-graph with $k\geq3$. Then the line graph of $G$ is symmetric with respect to graph entropy.
\end{thm3}
As a corollary to this result we show that the line graph of every bridgeless cubic graph is symmetric with respect to graph  entropy.

J. K\"{o}rner investigated the basic properties of the graph entropy in several papers from 1973 till 1992 (see J. K\"{o}rner \cite{JKor}-\cite{Jkor4}).

Let $F$ and $G$ be two graphs on the same vertex set $V$. Then the union of graphs $F$ and $G$ is the graph $F\cup G$ with vertex set $V$ and its edge set is the union of the edge set of graph $F$ and the edge set of graph $G$. That is
\begin{eqnarray}
&&V\left(F\cup G\right) = V,\nonumber\\
&&E\left(F\cup G\right) = E\left(F\right)\cup E\left(G\right).\nonumber
\end{eqnarray}
The most important property of the entropy of a graph is that it is sub-additive with respect to the union of graphs, that is
\[
H\left(F\cup G,P\right)\leq H\left(F,P\right) + H\left(G,P\right).
\]
This leads to the application of graph entropy for graph covering problem as well as the problem of perfect hashing.

The graph covering problem can be described as follows. Given a graph $G$ and a family of graphs $\mathcal G$ where each graph $G_i\in \mathcal G$ has the same vertex set as $G$, we want to cover the edge set of $G$ with the minimum number of graphs from $\mathcal G$. Using the sub-additivity of graph entropy one can obtain lower bounds on this number.

Graph entropy was used in a paper by Fredman and Koml\'{o}s for the minimum number of perfect hash functions of a given range that hash all $k$-element subsets of a set of a given size (see Fredman and Koml\'{o}s \cite{FK}).

\section{Preliminaries}\label{sec:basics}
\subsection{Entropy of Graphs}
Let $G$ be a graph on vertex set $V(G)=\{1,\cdots,n\}$, and $P=(p_1,\cdots,p_n)$ be a probability density on $V(G)$. The \emph{vertex packing polytope} of a graph $G$, i.e., $VP(G)$, is the convex hull of the characteristic vectors of its independent sets.

The \emph{entropy of $G$ with respect to $P$}, i.e., $H(G,P)$, is then defined as
\[
H(G,P) = \min_{\mathbf{a}\in VP(G)} \sum_{i=1}^n p_i\log (1/a_i).\label{eq:combent}
\]
\begin{remark}\label{rem:Rem1}
\emph{Note that the function $\sum_{i=1}^kp_i\log \frac{1}{a_i}$ in the definition of the graph entropy is a convex function and tends to infinity at the boundary of the non-negative orthant and tends monotonically to $-\infty$ along the rays from the origin.}\qed
\end{remark}
The main properties of graph entropy are \emph{monotonicity}, \emph{sub-additivity}, and \emph{additivity under vertex substitution}. Monotonicity is formulated in the following lemma (see G. Simonyi \cite{Simu}).

\begin{lemma}\emph{(J. K\"{o}rner).}\label{lem:mono}
Let $F$ be a spanning subgraph of a graph $G$. Then for any probability density $P$ we have $H(F,P)\leq H(G,P)$.
\end{lemma}

The sub-additivity was first recognized by K\"{o}rner in \cite{JKor1} and he proved the following lemma.
\begin{lemma}\label{lem:subadd}\emph{(J. K\"{o}rner).}
Let $F$ and $G$ be two graphs on the same vertex set $V$ and $F\cup G$ denote the graph on $V$ with edge set $E(F)\cup E(G)$. For any fixed probability density $P$ we have
\[
H\left(F\cup G,P\right) \leq H\left(F,P\right) + H\left(G,P\right).
\]
\end{lemma}
\begin{figure}[!t]%
\centering
\subfloat[A 5-cycle $G$.]{
\begin{tikzpicture}
[scale = 2]
            \draw (-.5,0) -- (-1,1);
            \draw (-1,1) -- (0,1.75);
            \draw (0,1.75) -- (1,1);
            \draw (1,1) -- (.5,0);
            \draw (.5,0) -- (-.5,0);
            \node[font=\small] at (0,2) {$u_1$};
            \node[font=\small] at (1.2,1.2) {$u_2$};
            \node[font=\small] at (.8,0) {$u_3$};
            \node[font=\small] at (-.8,0) {$u_4$};
            \node[font=\small] at (-1.3,1.3) {$u_5$};
            \filldraw [blue]
            (-0.5,0) circle (3pt)
            (-1,1) circle (3pt)
            (0,1.75) circle (3pt)
            (1,1) circle (3pt)
            (.5,0) circle (3pt);
\end{tikzpicture}
}
\qquad\qquad\qquad
\subfloat[A triangle $F$.]{
\begin{tikzpicture}
[scale = 2]
            \draw (-.5,1.5) -- (.5,1.5);
            \draw (.5,1.5) -- (0,2);
            \draw (0,2) -- (-.5,1.5);
            \node[font=\small] at (0,2.2) {$v_1$};
            \node[font=\small] at (.7,1.7) {$v_2$};
            \node[font=\small] at (-.7,1.7) {$v_3$};
            \filldraw [blue]
            (-0.5,1.5) circle (3pt)
            (.5,1.5) circle (3pt)
            (0,2) circle (3pt);
\end{tikzpicture}
}
\qquad\qquad\qquad
\subfloat[The graph $G_{u_1\longleftarrow F}$]{
\begin{tikzpicture}
[scale = 2]
            \draw (-.5,0) -- (-1,1);
            \draw (-1,1) to [out=120,in=150] (0,2);
            \draw (-1,1) -- (-.5,1.5);
            \draw (-1,1) -- (.5,1.5);
            \draw (0,2) to [out=30,in=60]  (1,1);
            \draw (-.5,1.5) -- (1,1);
            \draw (.5,1.5) -- (1,1);
            \draw (1,1) -- (.5,0);
            \draw (.5,0) -- (-.5,0);
            \draw (-.5,1.5) -- (.5,1.5);
            \draw (-.5,1.5) -- (0,2);
            \draw (.5,1.5) -- (0,2);
            \node[font=\small] at (0,2.2) {$v_1$};
            \node[font=\small] at (-.7,1.7) {$v_2$};
            \node[font=\small] at (.7,1.7) {$v_3$};
            \node[font=\small] at (1.2,1.2) {$u_2$};
            \node[font=\small] at (.8,0) {$u_3$};
            \node[font=\small] at (-.8,0) {$u_4$};
            \node[font=\small] at (-1.3,1.3) {$u_5$};
            \filldraw [blue]
            (-0.5,0) circle (3pt)
            (-1,1) circle (3pt)
            (0,2) circle (3pt)
            (-.5,1.5) circle (3pt)
            (.5,1.5) circle (3pt)
            (1,1) circle (3pt)
            (.5,0) circle (3pt);
\end{tikzpicture}
}
\caption{}%
\label{fig:Subs}%
\end{figure}
The notion of substitution is defined as follows. Let $F$ and $G$ be two vertex disjoint graphs and $v$ be a vertex of $G$. By substituting $F$ for $v$ we mean deleting $v$ and joining every vertex of $F$ to those vertices of $G$ which have been adjacent with $v$. We will denote the resulting graph $G_{v\leftarrow F}$. We extend this concept also to distributions. If we are given a probability distribution $P$ on $V(G)$ and a probability distribution $Q$ on $V(F)$ then by $P_{v\leftarrow Q}$ we denote the distribution on $V\left(G_{v\leftarrow F}\right)$ given by $P_{v\leftarrow Q}(x) = P(x)$ if $x \in V(G) \setminus {v}$ and $P_{v\leftarrow Q}(x) = P(x) Q(x)$ if $x \in V(F)$. This operation is illustrated in Figure \ref{fig:Subs}.

Now we state the following lemma whose proof can be found in J. K\"{o}rner, et. al. \cite{JKor2}.

\begin{lemma}\emph{(J. K\"{o}rner, G. Simonyi, and Zs.~Tuza).}\label{lem:Subs}
Let $F$ and $G$ be two vertex disjoint graphs, $v$ a vertex of $G$, while $P$ and $Q$ are probability distributions on $V(G) $ and $V(F)$, respectively. Then we have
\[
H\left(G_{v\leftarrow F}, P_{v\leftarrow Q}\right) = H\left(G,P\right) + P(v)H\left(F,Q\right).\qed
\]
\end{lemma}
Notice that the entropy of an empty graph (a graph with no edges) is always zero (regardless of the distribution on its vertices). Noting this fact, we have the following corollary as a consequence of Lemma \ref{lem:Subs}.

\begin{corollary}\label{cor:EntrDisc}
Let the connected components of the graph $G$ be the subgraphs $G_i$ and $P$ be a probability distribution on $V(G)$. If $x\in V(G_i)$, set
\[
P_i(x) = P(x)\left(P(V(G_i))\right)^{-1}, \quad x\in V(G_i).
\]
Then
\[
H\left(G,P\right) = \sum_i P\left(V(G_i)\right) H\left(G_i,P_i\right).
\]
\end{corollary}
Now we look at entropy of some graphs which are also mentioned in G. Simonyi \cite{Sim} and \cite{Simu} . The first one is the complete graph.

\begin{lemma}
For $K_n$, the complete graph on $n$ vertices, one has
\[
H\left(K_n,P\right) = H(P).\qed
\]
\end{lemma}
And the next one is the complete multipartite graph. Let $G = K_{m_1,m_2,\cdots,m_k}$ denote a complete $k$-partite graph with parts of size $m_1,m_2,\cdots,m_k$. Then we have the following lemma.

\begin{lemma}\label{lem:components}
Let $G = K_{m_1,m_2,\cdots,m_k}$. Given a distribution $P$ on $V(G)$ let $Q$ be the distribution on $\mathcal S(G)$, the set of maximal independent sets of $G$, given by $Q(J) = \sum_{x\in J}P(x)$ for each $J\in \mathcal S(G)$. Then $H(G,P) = H\left(K_k,Q\right)$.\qed
\end{lemma}
A special case of the above Lemma is the entropy of a complete bipartite graph with equal probability measure on its stable sets equal to 1.
Now, let $G$ be a bipartite graph with color classes $A$ and $B$. For a subset $D\subseteq A$, let $\mathcal{N}(D)$ denotes the the set of neighbors of $D$ in $B$, that is a subset of the vertices in $B$ which are adjacent to a vertex in $A$.

Given a distribution $P$ on $V(G)$ we have
\[
P(D) = \sum_{i\in D}p_i,\quad\forall D\subseteq V(G),
\]
Furthermore, defining the binary entropy as
\[
h(x) := -x\log x - (1 - x)\log (1 - x),\quad0\leq x\leq 1,
\]
J. K\"{o}rner and K. Marton proved the following theorem in \cite{JKor3}.

\begin{theorem}\emph{(J. K\"{o}rner and K. Marton).}\label{thm:bipentropy}
Let $G$ be a bipartite graph with no isolated vertices and $P$ be a probability distribution on its vertex set. If
\[
\frac{P(D)}{P(A)} \leq \frac{P(\mathcal{N}(D))}{P(B)},
\]
for all subsets $D$ of $A$, then
\[
H\left(G,P\right) = h\left(P(A)\right).
\]
And if
\[
\frac{P(D)}{P(A)} > \frac{P(\mathcal{N}(D))}{P(B)},
\]
then there exists a partition of $A = D_1\cup\cdots\cup D_k$ and a partition of $B = U_1\cup\cdots\cup U_k$ such that
\[
H\left(G,P\right) = \sum_{i=1}^k P\left(D_i \cup U_i\right) h\left(\frac{P(D_i)}{P(D_i\cup U_i)}\right).\qed
\]
\end{theorem}
\subsection{Minimum Entropy Colouring}\label{sec:chroment}
In this section, we explain minimum entropy colouring of the vertex set of a probabilistic graph $(G,P)$ which was previously studied by N. Alon and A. Orlitsky \cite{Alon96}.

Let $X$ be a random variable distributed over a countable set $V$ and $\pi$ be a partition of $V$, i.e., $\pi = \{C_1,\cdots,C_k\}$ and $V = \cup_{i = 1}^kC_i$. Then $\pi$ induces a probability distribution on its cells, that is

\[
p(C_i) = \sum_{v\in C_i} p(v), \quad\forall i\in\{1,\cdots,k\}.
\]

Therefore, the cells of $\pi$ have a well-defined entropy as follows:

\[
H\left(\pi\right) = \sum_{i = 1}^k p\left(C_i\right) \log \frac{1}{p\left(C_i\right)},
\]

If we consider $V$ as the vertex set of a probabilistic graph $(G,P)$ and $\pi$ as a partitioning of the vertices of $G$ into colour classes, then $H\left(\pi\right)$ is the entropy of a proper colouring of $V(G)$.

The \emph{chromatic entropy} of a probabilistic graph $(G,P)$ is defined as

\[
H_{\chi}(G,P) := \min\{H\left(\pi\right): \pi~\text{is a colouring of G}\},
\]

i.e. the lowest entropy of any colouring of $G$.
Consider a 5-cycle with two different probability distributions over its vertices, i.e., uniform distribution and another one given by
$p_1 = 0.3$, $p_2 = p_3 = p_5 = 0.2$, and $p_4 = 0.1$. In both of them we require three colours. In the first one, a colour is assigned to a single vertex and each of the other two colours are assigned to two vertices. Therefore, the chromatic entropy of the first probabilistic 5-cycle, i.e., $H(0.4,0.4,0.2)$ is approximately equal to 1.52.

For the second probabilistic 5-cycle, the chromatic entropy $
H(0.5,0.4,0.1)$ is approximately equal to 1.36. This chromatic entropy can be attained by choosing the colour classes to be $\{1,3\}$, $\{2, 5\}$, and $\{4\}$.

The following lemmas were proved in N. Alon and A. Orlitsky \cite{Alon96} and  J. Cardinal et.~al.~\cite{JC08}.

\begin{lemma}\emph{(N. Alon and A. Orlitsky).}\label{lem:H1}
Let $U$ be the uniform distribution over the vertices $V(G)$ of a probabilistic graph $(G,U)$ and $\alpha(G)$ be the independence number of the graph $G$. Then,
\[
H_{\chi}(G,U)\geq \log\frac{|V(G)|}{\alpha(G)}.\qed
\]
\end{lemma}

Let $\alpha(G,P)$ denote the maximum weight $P(S)$ of an independent set $S$ of a probabilistic graph $(G,P)$. Then we have the following lemma.
\begin{lemma}\emph{(J. Cardinal et.~al.).}\label{lem:H2}
For every probabilistic graph $(G,P)$, we have
\[
-\log\alpha(G,P)\leq H(G,P) \leq H_{\chi}(G,P) \leq \log\chi(G).\qed
\]
\end{lemma}

It may seem that non-uniform distribution decreases chromatic entropy $H_\chi(G,P)$, but the following example shows that this is not true.
Let us consider 7-star with $deg(v_1)=7$ and $deg(v_i) =1$ for $i\in\{2,\cdots,8\}$. If $p(v_1) = 0.5$ and $p(v_i) = \frac{1}{14}$ for $i\in\{2,\cdots,8\}$, then $H_\chi(G,P) = H(0.5,0.5) = 1$, while if $p(v_i) = \frac{1}{8}$ for $i\in\{1,\cdots,8\}$, then $H_\chi(G,P) = H(\frac{1}{8},\frac{7}{8})\leq H(0.5,0.5) = 1$.

Let $G_1,\cdots,G_n$ be graphs with vertex sets $V_1,\cdots,V_n$. The \emph{OR product}\index{product!OR product} of $G_1,\cdots,G_n$ is the graph $\bigvee_{i=1}^nG_i$ whose vertex set is $V^n$ and where two distinct vertices $(v_1,\cdots,v_n)$ and $(v^\prime_1,\cdots,v^\prime_n)$ are adjacent if for some $i\in\{1,\cdots,n\}$ such that $v_i\neq v^\prime_i$, $v_i$ is adjacent to $v^\prime_i$ in $G_i$. The $n$-fold OR product of $G$ with itself is denoted by $G^{\bigvee n}$.

N. Alon and A. Orlitsky \cite{Alon96} proved the following lemma which relates chromatic entropy to graph entropy.
\begin{lemma}\emph{(N. Alon and A. Orlitsky).}\label{lem:H0}
\[
\lim_{n\rightarrow\infty} \frac{1}{n} H_{\chi} (G^{\bigvee n}, P^{(n)}) = H(G,P).
\]
\end{lemma}

\subsection{Karush Kuhn Tucker (KKT) Conditions}\label{sec:kkt}
Karush Kuhn Tucker (KKT) optimality conditions in convex optimization are one of our tools in this paper. Thus, we explain these conditions briefly in this section. For a comprehensive explanation see Stephen Boyd, and Lieven Vanderberghe \cite{Boyd} section 5.5.3.
Let
\begin{eqnarray}
&&f_0(\mathbf x): \mathbb R^n\rightarrow\mathbb R,\nonumber\\
&&f_i(\mathbf x): \mathbb R^n\rightarrow\mathbb R,\quad i=1,\cdots,m,\nonumber\\
&&h_i(\mathbf x): \mathbb R^n\rightarrow\mathbb R,\quad i=1,\cdots,p.\nonumber
\end{eqnarray}
be convex differentiable functions.
We consider the following convex optimization problem,
\begin{eqnarray}
\text{minimize} &&f_0(\mathbf x)\label{eq:conv}\\
\text{subject to} && f_i(\mathbf x)\leq0,\quad i = 1,\cdots,m,\nonumber\\
&& h_i(\mathbf x) = 0,\quad i = 1,\cdots, p.\nonumber
\end{eqnarray}
Letting $\lambda_i$ and $\nu_i$ be the \emph{Lagrange multipliers} corresponding to constraints $f_i(\mathbf x)$ and $h_i(\mathbf x)$, respectively, we define the \emph{Lagrangian} $L:\mathbb R^n\times\mathbb R^m\times \mathbb R^p\rightarrow\mathbb R$ as
\[
L\left(\mathbf x,\lambda,\nu\right) = f_0(\mathbf x) + \sum_{i=1}^m\lambda_if_i(\mathbf x)+\sum_{i=1}^p\nu_ih_i(\mathbf x).
\]
Then $\mathbf x^*$ is an optimal solution to (\ref{eq:conv}) if and only if there exist Lagrange multipliers $\lambda^*$ and $\nu^*$ such that the following conditions hold,
\begin{eqnarray}
&&f_i(\mathbf x^*)\leq 0,\quad i=1,\cdots,m,\label{eq:kkt}\\
&&h_i(\mathbf x^*) = 0,\quad i =1,\cdots,p,\nonumber\\
&&\lambda_i^*\geq 0,\quad i=1,\cdots,m,\nonumber\\
&&\lambda_i^*f_i(\mathbf x^*) = 0,\quad i =1,\cdots,m\nonumber\\
&&\nabla f_0(\mathbf x^*) +\sum_{i=1}^m\lambda_i^*\nabla f_i(\mathbf x^*) +\sum_{i=1}^p\nu_i^*\nabla h_i(\mathbf x^*) = 0.\nonumber
\end{eqnarray}
The conditions (\ref{eq:kkt}) above are called \emph{Karush Kuhn Tucker (KKT)} optimality conditions.
\section{Graph Entropy and Fractional Chromatic Number}
In this section we investigate the relation between the entropy of a graph and its fractional chromatic number\index{fractional chromatic number} which was already established by G. Simonyi \cite{Simu}. First we recall that the \emph{fractional chromatic number} of a graph $G$ denoted by $\chi_f\left(G\right)$ is the minimum sum of nonnegative weights on the independent sets of $G$ such that for any vertex the sum of the weights on the independent sets of $G$ containing that vertex is at least one (see C. Godsil and G. Royle \cite{CGodsil} sections 7.1 to 7.5).
\begin{lemma}\emph{(G. Simonyi).}\label{lem:keylemma}
For a graph $G$ and probability density $P$ on its vertices with fractional chromatic number $\chi_f(G)$, we have
\[
\max_{P} H(G,P) = \log\chi_f(G).
\]
\end{lemma}
\proof
Note that for every graph $G$ we have $\left(\frac{1}{\chi_f(G)},\cdots,\frac{1}{\chi_f(G)}\right)\in VP(G)$. Thus for every probability density $P$, we have
\[
H\left(G,P\right)\leq \log \chi_f(G).
\]
Now, we show that graph $G$ has an induced subgraph $G^\prime$ with $\chi_f\left(G^\prime\right)=\chi_f\left(G\right)=\chi_f$ such that
if $\mathbf y\in VP\left(G^\prime\right)$ and $\mathbf y\geq\frac{\mathbf 1}{\chi_f}$, then $\mathbf y=\frac{\mathbf 1}{\chi_f}$.

Suppose the above statement does not hold for graph $G$. Consider all $\mathbf y\in VP(G)$ such that
\[
\mathbf y\geq\frac{\mathbf 1}{\chi_f(G)}.
\]
Note that there is not any $\mathbf y\in VP(G)$ such that
\[
\mathbf y>\frac{\mathbf 1}{\chi_f(G)},
\]
because then we have a fractional colouring with value strictly less than $\chi_f(G)$. Thus for every $\mathbf y\geq\frac{\mathbf 1}{\chi_f(G)}$ there is some $v\in V(G)$ such that $y_v=\frac{1}{\chi_f(G)}$.
For such a fixed $\mathbf y$, let
\[
\Omega_{\mathbf y} =\left\{v\in V(G):y_v>\frac{1}{\chi_f(G)}\right\}.
\]
Let $\mathbf y^*$ be one of those $\mathbf y$'s with $|\Omega_{\mathbf y}|$ of maximum size.
Let
\[
G^\prime = G\left[V(G)\setminus\Omega_{\mathbf y^*} \right].
\]
From our the definition of $G^\prime$ and fractional chromatic number, we have either
\[
\chi_f(G^\prime)<\chi_f(G).
\]
or
\[
\exists \mathbf y\in VP\left(G^\prime\right),\text{such that}~\mathbf y\geq\frac{\mathbf 1}{\chi_f}~\text{and}~\mathbf y\neq\frac{\mathbf 1}{\chi_f}.
\]
Suppose
\[
\chi_f(G^\prime)<\chi_f(G).
\]
Therefore
\[
\mathbf z = \frac{\mathbf{1}}{\chi_f(G^\prime)}\in VP(G^\prime)
\]
and consequently
\[
\mathbf z>\frac{\mathbf 1}{\chi_f(G)}.
\]
Without loss of generality assume that
\[
V(G)\setminus V(G^\prime) = \{1,\cdots, |V(G)\setminus V(G^\prime)|\}.
\]
Set
\begin{eqnarray}
&&\epsilon :=\frac{1}{2}\left(\min_{v\in\Omega_{\mathbf y^*}}y_v - \frac{1}{\chi_f(G)}\right)>0,\nonumber\\
&&\mathbf z^* = \left(\mathbf 0_{|V(G)\setminus V(G^\prime)|}^T,\mathbf z^T\right)^T\in VP(G).\nonumber
\end{eqnarray}
Then
\[
(1 - \epsilon)\mathbf y^*+\epsilon \mathbf z^*\in VP(G),
\]
which contradicts the maximality assumption of $\Omega_{\mathbf y^*}$. Thus, we have
\[
\chi_f(G^\prime) = \chi_f(G).
\]
Now we prove that if $y\in VP\left(G^\prime\right)$ and $\mathbf y\geq\frac{\mathbf 1}{\chi_f}$, then $\mathbf y=\frac{\mathbf 1}{\chi_f}$.

Suppose $\mathbf z^\prime$ be a point in $VP(G^\prime)$ such that $\mathbf z^\prime\geq \frac{1}{\chi_f}$ but $\mathbf z^\prime\neq \frac{1}{\chi_f}$.
Set
\[
\mathbf y^\prime = \left(\mathbf 0_{|V(G)\setminus V(G^\prime)|}^T,\mathbf z^{\prime^T}\right)^T\in VP(G).
\]
Then using the $\epsilon>0$ defined above, we have
\[
(1 - \epsilon)\mathbf y^*+\epsilon \mathbf y^\prime\in VP(G),
\]
which contradicts the maximality assumption of $\Omega_{\mathbf y^*}$.

Now, by I.Csisz\'{a}r et.~al.~\cite{Csis}, there exists a probability density $P^\prime$ on $VP(G^\prime)$ such that $H\left(G^\prime,P^\prime\right) = \log\chi_f$. Extending $P^\prime$ to a probability distribution $P$ as
\begin{equation}
p_i = \left\{ \begin{array}{rcl}
 p_i^\prime,& & i\in V(G) ,\\
0, && i\in V(G)\setminus V(G^\prime).
\end{array}\right.
\end{equation}
the lemma is proved. Indeed, suppose that $H(G,P)<H\left(G^\prime, P^\prime\right)$ and let $\bar{\mathbf y}\in VP(G)$ be a point in $VP(G)$ which gives $H\left(G,P\right)$. Let $\bar{\mathbf y}_{VP(G^\prime)}$ be the restriction of $\bar{\mathbf y}$ in $VP(G^\prime)$. Then there exists $\mathbf z\in VP(G^\prime)$ such that
\[
\mathbf z\geq\bar{\mathbf y}_{VP(G^\prime)}.
\]
This contradicts the fact that
\[
H\left(G^\prime, P^\prime\right) = \log \chi_f.
\]
\qed

\section{Symmetric Graphs}
A \emph{symmetric graph with respect to graph entropy} is a graph whose entropy is maximized with uniform probability distribution over its vertex set. In this section, we characterize different classes of graphs which are symmetric with respect to graph entropy. Particularly, we consider symmetric bipartite graphs, symmetric perfect graphs, and symmetric line graphs.

First, using chromatic entropy introduced in section \ref{sec:chroment}, we show that every vertex transitive graph is symmetric with respect to graph entropy.

\begin{theorem}\label{thm:vxtrans}
Let $G$ be a vertex transitive graph. Then the uniform distribution over vertices of $G$ maximizes $H\left(G,P\right)$.
That is $H\left(G,U\right) = \log \chi_f\left(G\right)$.
\end{theorem}

\proof
First note that for a vertex transitive graph $G$, we have $\chi_f\left(G\right) = \frac{|V(G)|}{\alpha(G)}$, and the $n$-fold OR product
$G^{\bigvee n}$ of a vertex transitive graph $G$ is also vertex transitive.  Now from Lemma \ref{lem:H1}, Lemma \ref{lem:H2} , and Lemma \ref{lem:keylemma},
we have

\begin{equation}
H\left(G^{\bigvee n},U\right) \leq \log \chi_f\left(G^{\bigvee n}\right) \leq H_{\chi}\left(G^{\bigvee n},U\right),\label{eq:U1}
\end{equation}

From N. Alon and A. Orltsky \cite{Alon96} and D. Ullman and E. Scheinerman \cite{fgt}, we have $H\left(G^{\bigvee n},U\right) = n H\left(G,U\right)$,
$\chi_f\left(G^{\bigvee n}\right) =  \chi_f\left(G\right)^n $, and $\log \chi_f\left(G\right)  = \lim_{n\rightarrow\infty} \frac{1}{n}\log\chi\left(G^{\bigvee n}\right)$. Hence, applying Lemma \ref{lem:H0} to equation (\ref{eq:U1}) and using squeezing theorem, we get

\begin{equation}
H\left(G,U\right) = \log \chi_f\left(G\right) = \lim_{n\rightarrow\infty} \frac{1}{n}\log\chi\left(G^{\bigvee n}\right) =
\lim_{n\rightarrow\infty} \frac{1}{n}H_\chi\left(G^{\bigvee n},U\right).\label{U2}
\end{equation}\qed

The following example shows that the converse of the above theorem is not true. Consider $G = C_4 \cup C_6$, with vertex sets $V(C_4) = \{v_1,v_2,v_3,v_4\}$ and $V(C_6) = \{v_5,v_6,v_7,v_8,v_9,v_{10}\}$, and parts $A = \{v_1,v_3,v_5,v_7,v_9\}$, $B = \{v_2,v_4,v_6,v_8,v_{10}\}$. Clearly, $G$ is not a vertex transitive graph, however, using Theorem \ref{thm:bipentropy}, one can see that the uniform distribution $U = \left(\frac{1}{10},\cdots,\frac{1}{10}\right)$ gives the maximum graph entropy which is $1$.

\begin{remark}
\emph{Note that the probability distribution which maximizes the graph entropy is not unique. Consider $C_4$ with vertex set $V(C_4) = \{v_1,v_2,v_3,v_4\}$ with parts $A = \{v_1,v_3\}$ and $B = \{v_2,v_4\}$. Using Theorem \ref{thm:bipentropy}, probability distributions $P_1 = (\frac{1}{4},\frac{1}{4},\frac{1}{4},\frac{1}{4})$ and $P_2 = (\frac{1}{8},\frac{1}{4},\frac{3}{8},\frac{1}{4})$ give the maximum graph entropy which is $1$.}
\end{remark}
\subsection{Symmetric Perfect Graphs}

Let $G=\left(V,E\right)$ be a graph. Recall that the fractional vertex packing polytope of $G$,i.e, $FVP(G)$ is defined as
\[
FVP(G) := \{\mathbf{x}\in\mathbb{R}_+^{|V|}: \sum_{v\in K}x_v\leq 1~\text{for all cliques K of G}\}.
\]
Note that for every graph $G$,  we have
\[
VP(G)\subseteq FVP(G).
\]
The following theorem was previously proved in V. Chv\'{a}tal \cite{chv} and D. R. Fulkerson \cite{Flker}.

\begin{theorem}
A graph $G$ is perfect if and only if $VP(G) = FVP(G)$.\qed
\end{theorem}

The following theorem which is called \emph{weak perfect graph theorem} is useful in the following discussion. This theorem was proved by  Lov\'{a}sz in \cite{Lov1} and \cite{Lov2} and is follows.
\begin{theorem}\label{thm:weakpf}
A graph $G$ is perfect if and only if its complement is perfect.\qed
\end{theorem}
Now, we prove the following theorem which is a generalization of our bipartite symmetric graphs with respect to graph entropy.
\begin{theorem}\label{thm:sympf}
Let $G=(V,E)$ be a perfect graph and $P$ be a probability distribution on $V(G)$. Then $G$ is symmetric with respect to graph entropy $H\left(G,P\right)$ if and only if $G$ can be covered by disjoint copies of its maximum-size cliques.
\end{theorem}

\proof

Suppose $G$ is covered by its maximum-sized cliques, say $Q_1,\cdots,Q_m$. That is $V(G) = V(Q_1) \dot{\cup}\cdots \dot{\cup}V(Q_m)$ and
$|V(Q_i)|=\omega(G),~\forall i\in [m]$.

Now, consider graph $T$ which is the disjoint union of the subgraphs induced by $V(Q_i)~\forall i\in [m]$. That $T = \dot{\bigcup}_{i=1}^m G\left[V(Q_i)\right]$. Noting that $T$ is a disconnected graph with $m$ components, using Corollary \ref{cor:EntrDisc} we have

\[
H\left(T,P\right) = \sum_i P(Q_i)H(Q_i,P_i).
\]

Now, having $V(T) = V(G)$ and $E(T)\subseteq E(G)$, we get $H\left(T,P\right)\leq H\left(G,P\right)$ for every distribution $P$. Using Lemma \ref{lem:keylemma}, this implies

\begin{equation}
H\left(T,P\right) = \sum_i P(Q_i)H\left(Q_i,P_i\right) \leq \log\chi_f(G),~\forall P,\label{eq:comps}
\end{equation}

Noting that $G$ is a perfect graph, the fact that complete graphs are symmetric with respect to graph entropy, $\chi_f\left(Q_i\right) = \chi_f(G) = \omega(G)=\chi(G),~\forall i\in [m]$, and (\ref{eq:comps}), we conclude that uniform distribution maximizes $H\left(G,P\right)$.

Now, suppose that $G$ is symmetric with respect to graph entropy. We prove that $G$ can be covered by its maximum-sized cliques. Suppose this is not true. We show that $G$ is not symmetric with respect to $H\left(G,P\right)$.

Denoting the minimum clique cover number of $G$ by $\gamma(G)$ and the maximum independent set number of $G$ by $\alpha(G)$, from perfection of $G$ and weak perfect theorem, we get $\gamma(G) = \alpha(G)$. Then, using this fact, our assumption implies that $G$ has an independent set $S$ with $|S|> \frac{|V(G)|}{\omega(G)}$.

We define a vector $\overline{\mathbf{x}}$ such that $\overline{x}_v = \frac{|S|}{|V|}$ if $v\in S$ and $\overline{x}_v = \frac{1 - \frac{|S|}{|V|}}{\omega - 1}$ if $v\in V(G) \backslash S$. Then, we can see that $\overline{\mathbf{x}}\in FVP(G) = VP(G)$. Let $t:=\frac{|S|}{|V|}$. Then, noting that $t>\frac{1}{\omega}$,

\begin{eqnarray}
H\left(G,U\right)&\leq& -\frac{1}{|V|}\sum_{v\in V} \log \overline{x}_v\nonumber\\
&=& -\frac{1}{|V|}\left(\sum_{v\in S}\log\overline{x}_v + \sum_{v\in V\backslash S}\overline{x}_v\right)\nonumber\\
&=&-\frac{1}{|V|}\left(|S|\log\alpha + (|V| - |S|)\log\frac{1-\alpha}{\omega-1}\right)\nonumber\\
&=&-t\log t - (1 - t)\log\frac{1 - t}{\omega - 1}\nonumber\\
&=&-t\log t - (\omega -1)\left(\frac{1-t}{\omega-1}\log\frac{1 - t}{\omega -1}\right) < \log \omega(G).\nonumber
\end{eqnarray}
\qed
Note that we have
\[
\gamma\left(G\right) = \alpha\left(G\right).
\]
The above theorem about symmetric perfect graphs is specialized to bipartite graphs in the following corollary. We give a separate proof for the following corollary, in S. S. C. Rezaei \cite{Saeed}, using Hall's theorem, K\"{o}nig's theorem (see D. West \cite{West} section 3.1) and Theorem \ref{thm:bipentropy}.
\begin{corollary}\emph{(Symmetric Bipartite Graphs).}\label{thm:Bip}
Let $G$ be a bipartite graph with parts $A$ and $B$, and no isolated vertices. Then, uniform probability distribution $U$ over the vertices of $G$ maximizes $H\left(G,P\right)$ if and only if $G$ has a perfect matching.
\end{corollary}
We also have the following corollary.
\begin{corollary}
Let $G$ be a connected regular line graph without any isolated vertices with valency $k>3$. Then if $G$ is covered by copies of its disjoint maximum-size cliques, then $G$ is symmetric with respect to $H(G,P)$.
\end{corollary}

\proof

Let $G = L(H)$ for some graph $H$. Then either $H$ is bipartite or regular. If $H$ is bipartite, then $G$ is perfect (see D. West \cite{West} sections 7.1 and 8.1) and because of Theorem \ref{thm:sympf} we are done. So suppose that $H$ is not bipartite. Then each clique of size $k$ in $G$ corresponds to a vertex $v$ in $V(H)$ and the edges incident to $v$ in $H$ and vice versa. That is because any such cliques in $G$ contains a triangle and there is only one way extending that triangle to the whole clique which corresponds to edges incident with the corresponding vertex in $H$. This implies that the covering cliques in $G$ give an independent set in $H$ which is also a vertex cover in $H$. Hence $H$ is a bipartite graph and hence $G $ is perfect. Then due to Theorem \ref{thm:sympf} the theorem is proved.
\qed
Now, considering that finding the clique number of a perfect graph can be computed in polynomial time and using the weak perfect graph theorem, we conclude that one can decide in polynomial time whether a perfect graph is symmetric with respect to graph entropy.
\subsection{Symmetric Line Graphs}

Let $G_2$ be a line graph of some graph $G_1$, i.e, $G_2 = L(G_1)$.  Let $|V(G_1)|=n$ and $|E(G_1)| = m$.
Note that every matching in $G_1$ corresponds to an independent set in $G_2$ and every independent set in $G_2$ corresponds to a matching in $G_1$. Furthermore, the fractional edge-colouring number of $G_1$, i.e., $\chi_f^\prime(G_1)$ is equal to the fractional chromatic number of $G_2$, i.e.,$\chi_f(G_2)$. Thus
\[
\chi_f^\prime(G_1) = \chi_f(G_2).
\]
Moreover, note that the vertex packing polytope $VP(G_2)$ of $G_2$ is the matching polytope $MP(G_1)$ of $G_1$ (see L. Lov\'{a}sz and M. D. Plummer \cite{Lov3} chapter 12). That is
\[
VP(G_2) = MP(G_1).
\]

These facts motivate us to study line graphs which are symmetric with respect to graph entropy.

We recall that a vector $\mathbf x\in \mathbb R_+^m$ is in the matching polytope $MP(G_1)$ of the graph $G_1$ if and only if it satisfies (see A. Schrijver \cite{Schriv1}).
\begin{eqnarray}\label{eq:MatchPoly}
&x_e\geq 0~~~~~~~~ &\forall e\in E(G_1),\nonumber\\
&x(\delta(v))\leq 1 ~~~~~~~~~&\forall v\in V(G_1),\\
&x\left(E[U]\right)\leq\lfloor\frac{1}{2}|U|\rfloor, ~~~~~~~~~&\forall U\subseteq V(G_1)~\text{with}~|U|~\text{odd.}\nonumber
\end{eqnarray}
Let $\mathcal M$ denote the family of all matchings in $G_1$, and for every matching $M\in\mathcal M$ let the charactersitic vector $\mathbf b_M\in\mathbb R_+^m$ be as
\begin{equation}
(\mathbf b_M)_e = \left\{ \begin{array}{rcl}
1, & & e\in M ,\\
0, && e\notin M.
\end{array}\right.
\end{equation}
Then the \emph{fractional edge-colouring number}\index{fractional edge-colouring number} $\chi_f^\prime(G_1)$ of $G_1$ is defined as
\[
\chi_f^\prime(G_1):=\min\{\sum_{M\in\mathcal M}\lambda_M|\mathbf{\lambda}\in\mathbb R_+^{\mathcal M},~\sum_{M\in\mathcal M}\lambda_M\mathbf b_M = \mathbf 1\}.
\]
If we restrict $\lambda_M$ to be an integer, then the above definition give rise to the edge colouring number of $G_1$, i.e., $\chi^\prime(G_1)$. Thus
\[
\chi_f^\prime(G)\leq \chi^\prime(G).
\]
As an example considering $G_1$ to be the Petersen graph, we have
\[
\chi_f^\prime(G)=\chi^\prime(G) = 3.
\]

The following theorem which was proved by Edmonds, gives a characterization of the fractional edge-colouring number $\chi_f^\prime(G_1)$ of a graph $G_1$ (see A. Schrijver \cite{Schriv1} section 28.6).
\begin{theorem}\label{thm:edgecol1}
Let $\Delta(G_1)$ denote the maximum degree of $G_1$. Then the fractional edge-colouring number of $G_1$ is obtained as
\[
\chi_f^\prime(G_1) = \max\{\Delta(G_1),\max_{U\subseteq V,~|U|\geq3}\frac{|E(U)|}{\lfloor\frac{1}{2}|U|\rfloor}\}.\qed
\]
\end{theorem}
Following A. Schrijver \cite{Schriv1} we call a graph $G_1$ a \emph{$k$-graph}\index{$k$-graph} if it is $k$-regular and its fractional edge coloring number $\chi_f^\prime(H)$ is equal to $k$.
The following corollary characterizes a $k$-graph (see A. Schrijver \cite{Schriv1} section 28.6).
\begin{corollary}\label{cor:kgraph}
Let $G_1=(V_1,E_1)$ be a $k$-regular graph. Then $\chi^\prime_f(G_1)=k$ if and only if $|\delta(U)|\geq k$ for each odd subset $U$ of $V_1$.
\end{corollary}
The following theorem introduces a class of symmetric line graphs with respect to graph entropy. The main tool in the proof of the following theorem is Karush-Kuhn-Tucker (KKT) optimality conditions in convex optimization explained in section \ref{sec:kkt}.
\begin{theorem}
Let $G_1$ be a $k$-graph with $k\geq3$. Then the line graph $G_2 = L(G_1)$ is symmetric with respect to graph entropy.
\end{theorem}
\proof

From our discussion above we have
\[
H\left(G_2,P\right) = \min_{\mathbf x\in MP(G_1) }\sum_{e\in E(G_1)} p_e\log \frac{1}{x_e},
\]
Let $\lambda_v,~\gamma_U\geq 0$ be the Lagrange multipliers corresponding to inequalities $x(\delta(v))\leq 1$ and $x\left(E[U]\right)\leq\lfloor\frac{1}{2}|U|\rfloor$ in the description of the matching polytope $MP(G_1)$ in (\ref{eq:MatchPoly}) for all $v\in V(G_1)$ and for all $U\subseteq V(G_1)$ with $|U|$ odd, and $|U|\geq 3$, repectively. From our discussion in Remark \ref{rem:Rem1}, the Lagrange mulitipliers corresponding to inequalities $x_e\geq 0$ are all zero.

Set
\[
g(\mathbf x) = -\sum_{e\in E(G_1)} p_e\log x_e,
\]
Then the Lagrangian of $g(\mathbf x)$ is
\begin{eqnarray}
L\left(\mathbf x, \mathbf{\lambda}, \mathbf{\gamma}\right) &&=-\sum_{e\in E(G_1)}p_e\log x_e +\sum_{e=\{u,v\}}\left(\lambda_u+\lambda_v\right)\left(x_e - 1\right)\nonumber\\
&&+\sum_{e\in E(G_1)}\sum_{\substack{U\subseteq V,\\U\ni e,|U|~\text{odd},~|U|\geq 3 }}\gamma_U x_e - \sum_{\substack{U\subseteq V,\\|U|~\text{odd},~|U|\geq 3}}\lfloor\frac{1}{2}|U|\rfloor,
\end{eqnarray}
Using KKT conditions (see S. Boyd, and L. Vanderberghe \cite{Boyd} section 5.5.3), the vector $\mathbf x^*$ minimizes $g(\mathbf x)$ if and only if it satisfies \begin{eqnarray}\label{eq:KKT2}
&&\frac{\partial L}{\partial x_e^*} = 0,\nonumber\\
&&\rightarrow -\frac{p_e}{x_e^*} + \left(\lambda_u+\lambda_v\right) + \sum_{\substack{U\subseteq V,\\U\ni e, |U|~\text{odd},~|U|\geq 3}}\gamma_U = 0~\text{for}~e=\{u,v\}.
\end{eqnarray}
Fix the probability density to be uniform over the edges of $G_1$, that is
\[
p_e = \frac{1}{m}, ~~~\forall e\in E(G_1),
\]
Note that the vector $\frac{\mathbf 1}{k}$ is a feasible point in the matching polytope $MP(G_1)$. Now, one can verify that specializing the variables as
\begin{eqnarray}
&&\mathbf x^* = \frac{\mathbf 1}{k},\nonumber\\
&&\gamma_U = 0~~~~~~~~~~~\forall U\subseteq V,~|U|~\text{odd},~|U|\geq3\nonumber\\
&&\lambda_u=\lambda_v = \frac{k}{2m}~~~~~~~~~~\forall~e=\{u,v\}.\nonumber
\end{eqnarray}
satisfies the equations (\ref{eq:KKT2}). Thus
\[
H\left(G_2,U\right) = \log k.
\]
Then using Lemma \ref{lem:keylemma} and the assumption $\chi_f(G_2) = k$ the theorem is proved.
\qed

It is well known that cubic graphs has a lot of interesting structure. For example, it can be checked that every edge in a bridgeless cubic graph is in a perfect matching.
Now we have the following interesting statement for every cubic bridgeless graph.
\begin{corollary}\label{cor:symmcubic1}
The line graph of every cubic bridgeless graph $G_1=(V_1,E_1)$ is symmetric with respect to graph entropy.
\end{corollary}
\proof
We may assume that $G_1$ is connected. Let $U\subseteq V_1$ and let $U_1\subseteq U$ consist of vertices $v$ such that $\delta(v)\cap\delta(U)=\emptyset$. Then using handshaking lemma for $G_1[U]$, we have
\[
3|U_1| + 3|U\setminus U_1| - |\delta(U)| = 2|E(G_1[U])|.
\]
And consequently,
\[
3|U| = |\delta(U)| \mod 2,
\]
Assuming $|U|$ is odd and noting that $G_1$ is bridgeless, we have
\[
\delta(U)\geq 3.
\]
Then, considering Corollary \ref{cor:kgraph}, the corollary is proved.
\qed

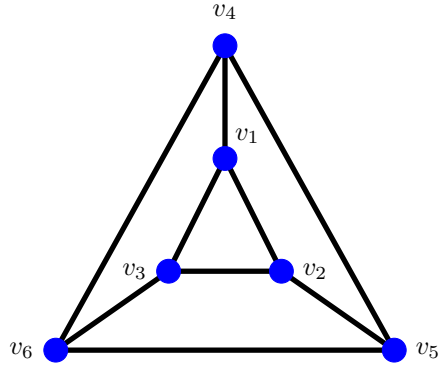
\begin{figure}[!t]
    \begin {center}
        \begin{tikzpicture}
        [scale = 1.5,
        foo/.style={line width = 15pt}]
            \draw[line width=2pt] (0,1) -- (0.5,0);
            \draw[line width=2pt] (0.5,0) -- (-.5,0);
            \draw[line width=2pt] (-.5,0) -- (0,1);
            \draw[line width=2pt] (0,1) -- (0,2);
            \draw[line width=2pt] (0,2) -- (-1.5,-.7);
            \draw[line width=2pt] (-1.5,-.7) -- (1.5,-.7);
            \draw[line width=2pt] (1.5,-.7) -- (0,2);
            \draw[line width=2pt] (1.5,-.7) -- (.5,0);
            \draw[line width=2pt] (-1.5,-.7) -- (-.5,0);
            \node[font=\small] at (0.2,1.2) {$v_1$};
            \node[font=\small] at (.8,0) {$v_2$};
            \node[font=\small] at (-.8,0) {$v_3$};
            \node[font=\small] at (0,2.3) {$v_4$};
            \node[font=\small] at (1.8,-.7) {$v_5$};
            \node[font=\small] at (-1.8,-.7) {$v_6$};
            \filldraw [blue]
            (0,1) circle (3pt)
            (.5,0) circle (3pt)
            (-.5,0) circle (3pt);
            \filldraw [blue]
            (0,2) circle (3pt)
            (1.5,-.7) circle (3pt)
            (-1.5,-.7) circle (3pt);
        \end{tikzpicture}
    \end{center}
\caption{A bridgeless cubic graph.}\label{fig:bridgelesscubic}
  \end{figure}
Figure \ref{fig:bridgelesscubic} shows a bridgeless cubic graph which is not edge transitive and its edges are not covered by disjoint copies of stars and triangles. Thus the line graph of the shown graph in Figure \ref{fig:bridgelesscubic} is neither vertex transitive nor covered by disjoint copies of its maximum size cliques. However, it is symmetric with respect to graph entropy by Corollary \ref{cor:symmcubic1}.

Figure \ref{fig:cubicgraph} shows a cubic graph with a bridge. The fractional edge chromatic number of this graph is $3.5$ while the entropy of its line graph is $1.75712$, i.e., $\log_23.5 = 1.8074>1.75712$. Thus, its line graph is not symmetric with respect to graph entropy, and we conclude that Corollary \ref{cor:symmcubic1} is not true for cubic graphs with bridges.
\begin{figure}[!t]\label{fig:cubicgraph}
    \begin {center}
        \begin{tikzpicture}
        [scale = 1.1,
        foo/.style={line width = 15pt}]
            \draw[line width=2pt] (-5,0) -- (-3,0);
            \draw[line width=2pt] (-5,0) -- (-2,1.5);
            \draw[line width=2pt] (-5,0) -- (-2,-1.5);
            \draw[line width=2pt] (-3,0) -- (-2,1.5);
            \draw[line width=2pt] (-3,0) -- (-2,-1.5);
            \draw[line width=2pt] (-2,1.5) -- (-1,0);
            \draw[line width=2pt] (-2,-1.5) -- (-1,0);
            \draw[line width=2pt] (-1,0) -- (1,0);
            \draw[line width=2pt] (1,0) -- (2,1.5);
            \draw[line width=2pt] (1,0) -- (2,-1.5);
            \draw[line width=2pt] (2,1.5) -- (3,0);
            \draw[line width=2pt] (2,-1.5) -- (3,0);
            \draw[line width=2pt] (3,0) -- (5,0);
            \draw[line width=2pt] (2,1.5) -- (5,0);
            \draw[line width=2pt] (2,-1.5) -- (5,0);

            \node[font=\small] at (-5.3,0) {$v_1$};
            \node[font=\small] at (-2,1.8) {$v_2$};
	    \node[font=\small] at (-2.7,0) {$v_3$};
           \node[font=\small] at (-2,-1.8) {$v_4$};
           \node[font=\small] at (-.7,.3) {$v_5$};
           \node[font=\small] at (.7,.3) {$v_6$};
           \node[font=\small] at (2,1.8) {$v_7$};
           \node[font=\small] at (2,-1.8) {$v_8$};
           \node[font=\small] at (2.7,0) {$v_9$};
           \node[font=\small] at (5.4,0) {$v_{10}$};

            \filldraw [blue]
            (-5,0) circle (3pt)
	    (-2,1.5) circle (3pt)
            (-3,0) circle (3pt)
            (-2,-1.5) circle (3pt)
            (-1,0) circle (3pt)
            (1,0) circle (3pt)
            (2,1.5) circle (3pt)
            (2,-1.5) circle (3pt)
            (3,0) circle (3pt)
            (5,0) circle (3pt);

        \end{tikzpicture}
    \end{center}
\caption{A cubic one-edge connected graph.}\label{fig:cubicgraph}
  \end{figure}
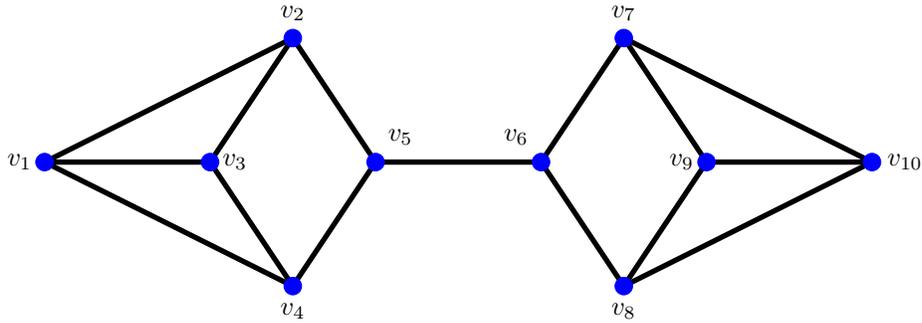
\section{Open Questions}
In this section we point out some of the problems that are worth considering for future research.
\subsection{Symmetric Graphs}
We defined symmetric graphs with respect to graph entropy in this paper. As the main result of this paper, we characterized symmetric perfect graphs. Furthermore, we proved that there are some other classes of graphs such as vertex transitive graphs and line graph of bridgeless cubic graphs which are symmetric with respect to graph entropy. From what discussed above and the symmetric graphs considered in this paper, the next natural class of graphs to consider is \emph{strongly regular graphs}. It is worth investigating if strongly regular graphs are symmetric with respect to graph entropy. For the study of the structural properties of strongly regular graphs see Godsil and Royle\cite{CGodsil} chapter 10. Furthermore, it is worth noting that finding the main properties of symmetric graphs with respect to graph entropy is another interesting open problem.
\subsection{Normal Graph Conjecture}
Let $G$ be a graph. A set $\mathcal A$ of subsets of $V(G)$ is a \emph{covering}, if every vertex of $G$ is contained in an element of $\mathcal A$.

We say that graph $G$ is \emph{Normal} if there exists two coverings $\mathcal C$ and $\mathcal S$ such that every element $C$ of $\mathcal C$ is a clique and every element $S$ of $\mathcal S$ is an independent set and the intersection of any element of $\mathcal C$ and any element of $\mathcal S$ is nonempty, i.e.,
\[
C\cap S \neq \emptyset,~\forall C\in \mathcal C,~S\in\mathcal S.
\]
Recall from the sub-additivity of Graph Entropy, we have
\begin{equation}\label{eq:subadd6}
H(P)\leq H(G,P) + H(\overline G,P).
\end{equation}
A probabilistic graph $(G,P)$ is \emph{weakly splitting} if there exists a nowhere zero probability distribution $P$ on its vertex set which makes inequality (\ref{eq:subadd6}) equality. The following lemma was proved in J. K\"{o}rner et. al. \cite{JKor2}.
\begin{lemma}\label{lem:normal1}\emph{(J. K\"{o}rner, G. Simonyi, and Zs. Tuza)}
A graph $G$ is weakly splitting if and only if it is normal.
\end{lemma}
It is known that every perfect graph is also a normal graph (see J. K\"{o}rner \cite{JKor01}). The following conjecture was proposed in C. De Simone and J. K\"{o}rner \cite{De}.
\begin{conjecture}\emph{(Normal Graph Conjecture).}
A graph is hereditarily normal if and only if the graph nor its complement contains $C_5$ or $C_7$ as an induced subgraph.
\end{conjecture}
A \emph{circulant} $C_n^k$ is a graph with vertex set $\{1,\cdots,n\}$, and two vertices $i\neq j$ are adjacent if and only if
\[
i - j \equiv k~\mathrm{mod}~n.
\]
We assume $k\geq 1$ and $n\geq 2(k+1)$ to avoid cases where $C_n^k$ is an independent set or a clique. A. K. Wagler \cite{Wag} proved the Normal Graph Conjecture for circulants $C_n^k$.

One direction for future research is investigating the Normal Graph Conjecture for general circulants and Cayley graphs.
\subsection{Graph Entropy and Graph Homomorphism}
Given a probabilistic graph $(G,P)$, K. Marton in K. Marton \cite{Mart} introduced a functional $\lambda(G,P)$ which is analogous to Lov\'{a}sz's bound $\vartheta(G)$ on Shannon capacity of graphs. Similar to $\vartheta(G)$, the probabilistic functional $\lambda(G,P)$ is based on the concept of \emph{orthonormal representation} of a graph which is recalled here.

Let $U=\{\mathbf u_i:i\in V(G)\}$ be a set of unit vectors of a common dimension $d$ such that
\[
\mathbf u_i^T\mathbf u_j = 0~\text{if}~i\neq j~\text{and}~\{i,j\}\notin E(G).
\]
Let $\mathbf c$ be a unit vector of dimension $d$. Then, the system $\left(U,\mathbf c\right)$ is called an orthonormal representation of the graph $G$ with handle $\mathbf c$.

Letting $T(G)$ denote the set of all orthonormal representations for graph $G$ with a handle $\mathbf c$, Lov\'{a}sz \cite{Lov} defined
\[
\vartheta(G)=\min_{\left(U,\mathbf c\right)\in T(G)}\max_{i\in V(G)}\frac{1}{(\mathbf u_i,\mathbf c)^2}.
\]
Then it is shown in Lov\'{a}sz \cite{Lov} that zero error Shannon capacity $C(G)$ can be bounded above by $\vartheta(G)$ as
\[
C(G)\leq\log \vartheta(G).
\]
Let $P$ denote the probability distribution over the vertices of $G$, and $\epsilon>0$.
Then the \emph{capacity of the graph relative to} $P$ is
\[
C(G,P)=\lim_{\epsilon\rightarrow0}\limsup_{n\rightarrow\infty}\frac{1}{n}\log
\alpha\left(G^(P,\epsilon)\right).
\]
A probabilistic version of $\vartheta(G)$ denoted by $\lambda(G,P)$ is defined in K. Marton \cite{Mart} as
\[
\lambda(G,P):=\min_{\left(U,\mathbf c\right)\in T(G)}\sum_{i\in V(G)}P_i\log\frac{1}{(\mathbf u_i,\mathbf c)^2}.
\]
K. Marton \cite{Mart} showed that
\begin{theorem}\emph{(K. Marton)}
The capacity of a probabilistic graph $(G,P)$ is bounded above by $\lambda(G,P)$, i.e.,
\[
C(G,P)\leq \lambda(G,P).
\]
\end{theorem}
The following theorem was proved in K. Marton \cite{Mart} which relates $\lambda(G,P)$ to $H(G,P)$.
\begin{theorem}\emph{(K. Marton)}
For any probabilistic graph $(G,P)$,
\[
\lambda\left(\overline G, P\right)\leq H\left(G,P\right).
\]
Furthermore, equality holds if and only if $G$ is perfect.
\end{theorem}
K. Marton \cite{Mart} also related $\lambda(G,P)$ to $\vartheta(G)$ by showing
\begin{equation}\label{eq:Hom1}
\max_{P}\lambda(G,P) = \log\vartheta(G).
\end{equation}
It is worth mentioning that $\vartheta(G)$ can be defined in terms of graph homomorphisms as follows.

Let $d\in\mathbb{N}$ and $\alpha<0$. Then we define $S(d,\alpha)$ to be an infinite graph whose vertices are unit vectors in $\mathbb{R}^d$. Two vertices $\mathbf u$ and $\mathbf v$ are adjacent if and only if $\mathbf u\mathbf v^T = \alpha$. Then
\begin{equation}\label{eq:Hom2}
\vartheta\left(\overline G\right)=\min\left\{1 - \frac{1}{\alpha}:G\rightarrow S\left(d,\alpha\right),~\alpha<0\right\}.
\end{equation}
Thus, noting (\ref{eq:Hom1}) and (\ref{eq:Hom2}) and the above discussion, investigating the relationship between graph homomorphism and graph entropy which may lead to investigating the relationship between graph homomorphism and graph covering problem seems interesting.

 \end{document}